\documentclass[11pt]{amsart}

\usepackage{amsmath}
\usepackage{amssymb}
\usepackage{graphicx}
\usepackage{amsthm}
\usepackage{enumerate}
\usepackage{mathrsfs}
\usepackage{upgreek}
\usepackage[mathscr]{eucal}
\newtheorem{theorem}{Theorem}[section]
\newtheorem*{thm-non}{Theorem}
\newtheorem{lemma}[theorem]{Lemma}

\newtheorem*{exp-non}{Examples}
\newtheorem{notation}[theorem]{Notation}

\DeclareMathOperator*{\argmax}{argmax}
\DeclareMathOperator*{\argmin}{argmin}

\theoremstyle{plain}
\newcommand{\thistheoremname}{}
\newtheorem{genericthm}[theorem]{\thistheoremname}
\newenvironment{namedtheorem}[1]
{\renewcommand{\thistheoremname}{#1}%
	\begin{genericthm}}
	{\end{genericthm}}

\numberwithin{equation}{section}
\title[The novel Tauberian conditions associated with the $(\overline{N},p,q)$ summability]
{The novel Tauberian conditions associated with the $(\overline{N},p,q)$ summability of double sequences}

\author[Z. \"{O}nder]{Zerrin \"{O}nder}
\address[Z. \"{O}nder]{U\c{s}ak University \\ Department of Mathematics \\ 64000 U\c{s}ak Turkey}
\email{{\tt zerrin.onder11@gmail.com}}

\author[E. Sava\c{s}]{Ekrem Sava\c{s}}
\address[E. Sava\c{s}]{U\c{s}ak University \\ Department of Mathematics \\ 64000 U\c{s}ak Turkey}
\email{\tt ekrem.savas@usak.edu.tr}

\author[\.{I}. \c{C}anak]{\.{I}brahim \c{C}anak}
\address[\.{I}. \c{C}anak]{Ege University \\ Department of Mathematics \\ 35100 \.{I}zmir Turkey}
 \email{\tt ibrahim.canak@ege.edu.tr}

\keywords{Double sequences, convergence in Pringsheim's sense, $(\overline{N},p,q)$ summability, regularly varying sequences, slowly decreasing sequences, slowly oscillating sequences, Tauberian conditions and theorems, weighted mean summability method}

\subjclass[2010]{40A05, 40E05, 40G99}

\dedicatory{Dedicated to Ravi P. Agarwal on his 76th birthday}
\begin{document}

	\begin{abstract}
	In this paper, our primary objective is to provide a fresh perspective on the relationship between the $(\overline{N},p,q)$ method, which is a product of relevant one-dimensional summability methods, and $P$-convergence for double sequences. To accomplish this objective, we establish certain Tauberian conditions that control the behavior of a double sequence in terms of both $O_L$-oscillation and $O$-oscillation in certain senses, building a bridge between $(\overline{N},p,q)$ summability and $P$-convergence, while imposing certain restrictions on the weight sequences. As special circumstances of our findings, we demonstrate that Landau-type $O_L$ condition with respect to $(P_m)$ and $(Q_n),$ as well as Hardy-type $O$ condition with respect to $(P_m)$ and $(Q_n),$ serve as Tauberian conditions for $(\overline{N},p,q)$ summability under particular additional conditions. Consequently, these results encompass all classical Tauberian theorems, including conditions such as slow decrease or slow oscillation in certain senses.
\end{abstract}

\maketitle
%%%%%%%%%%%%%%%%%%%%%%%%%%%%%%%%%%%%%%%%%

\section{Introduction}

In the early 1900s, when the concept of convergence in Pringsheim's sense, known as $P$-convergence, came to light, the evolvement of summability theory for single sequences to multiple ones was still in its infancy. The detailed exploration of double sequences by Hardy \cite{hardy} and Bromwich \cite{bromwich} marked a significant turning point, leading to a surge of interest in this novel type of sequences. The application of the weighted mean methods to double sequences began, to the best of our knowledge, with the pioneering work of Baron and Stadtm\"{u}ller \cite{baronstadtmuller}. In their seminal paper \cite{baronstadtmuller}, Baron and Stadtm\"{u}ller embarked on an analysis of the relationship between the $(\overline{N},p,q)$ method, which is a product of relevant one-dimensional methods, and $P$-convergence for double sequences. They established that the necessary conditions for (boundedly) $P$-convergence of a double sequence, which is (boundedly) $(\overline{N},p,q)$ summable, are Hardy type $O$-condition relative to $(P_m)$ and $(Q_n),$ where the weights $(P_m)$ and $(Q_n)$ are regularly varying.
%\begin{equation*}
	%\sup_{n\in\mathbb{N}}\left(\Delta_{10}u_{mn}\right)=O\left(\frac{p_m}{P_m}\right)\,\,\,\,\,\text{and}\,\,\,\,\,
	%\sup_{m\in\mathbb{N}}\left(\Delta_{01}u_{mn}\right)=O\left(\frac{q_n}{Q_n}\right)
	%\end{equation*}
Instead of using factorable weights, as was done previously, Stadtm\"{u}ller \cite{stadtmuller2} extended these conditions to non-factorable weights, and so, both generalized $O_L$-Tauberian conditions given by M\'{o}ricz \cite{moricz1994} for $(C, 1, 1)$ method and demonstrated that the conditions could be relaxed. Following Stadtmüller's work \cite{stadtmuller2}, Chen and Hsu \cite{chenhsu} established Tauberian theorems for double sequences, dealing with implication from $(\overline{N},p,q)$ summability to $P$-convergence under Landau-type conditions, Schmidt-type slow decrease conditions and more general conditions involving the concept of deferred means. Reducing the assumptions asserted by Stadtm\"{u}ller \cite{stadtmuller2}, M\'{o}ricz and Stadtm\"{u}ller \cite{moriczstadtmuller} further examined the conditions necessary for (boundedly) $(\overline{N},p,q)$ summable double sequences to be (boundedly) $P$-convergent by utilizing the classes $\Lambda_u$ and $\Lambda_{\ell},$ constructed based on non-factorable weights. Lastly, Belen \cite{belen} introduced the concept of double weighted generator sequences in certain senses, which represent the difference between double sequences and their $(\overline{N},p,q)$ means. Belen \cite{belen} also identified certain conditions formed via these sequences, such as $P_{m-1}\Delta_{10}V_{mn}^{11^{(0)}}\left(\Delta_{11}(u)\right)=O_L(p_m)$ and $Q_{n-1}\Delta_{01}V_{mn}^{11^{(0)}}\left(\Delta_{11}(u)\right)=O_L(q_n)$, as Tauberian conditions for the $(\overline{N},p,q)$ method, subject to additional conditions imposed on the weights $(P_m)$ and $(Q_n).$
%\begin{equation*}
	%\Delta_{10}V_{mn}^{11^{(0)}}\left(\Delta_{11}(u)\right)=O_L\left(\frac{p_m}{P_{m-1}}\right)\,\,\,\,\,\text{and}\,\,\,\,\,
	%\Delta_{01}V_{mn}^{11^{(0)}}\left(\Delta_{11}(u)\right)=O_L\left(\frac{q_n}{Q_{n-1}}\right)
	%\end{equation*}

%%%%%%%%%buradan devam%%%%%%%%%%%%

Together with the mentioned studies so far, the matter that urges us to make this work is the idea of furthering the results obtained by Boos \cite{boos} for single sequences by extending to $(\overline{N},p,q)$ summable double sequences. In \cite{boos}, Boos formulated these results as follows:
\begin{theorem} \label{A}
	Let $(p_n)$ be a sequence which has the property
	$$\frac{P_n}{P_{n+1}}\to 1\,\,\,\,\text{as}\,\,\,\,n\to\infty.$$
	If a sequence $(u_{n})$ of real numbers is $(\overline{N},p)$ summable to $\ell$ and slowly decreasing relative to $(P_n),$ then  $(u_{n})$ is convergent to $\ell$.
\end{theorem}
\begin{theorem} \label{B}
	Let $(p_n)$ be a sequence which has the property
	$$\frac{P_n}{P_{n+1}}\to 1\,\,\,\,\text{as}\,\,\,\,n\to\infty.$$
	If a sequence $(u_{n})$ of complex numbers is $(\overline{N},p)$ summable to $\ell$ and slowly oscillating relative to $(P_n),$ then  $(u_{n})$ is convergent to $\ell$.
\end{theorem}
\noindent One of purposes of this paper is to extend Theorems \ref{A} and \ref{B} given for $(\overline{N},p)$ summable sequences of real and complex numbers to $(\overline{N},p,q)$ summable double sequences of real and complex numbers. The other is to indicate that our results obtained in this paper include all of the classical Tauberian theorems for double sequences which $P$-convergence follows from Ces\`{a}ro and logarithmic summability under slow decrease (or oscillation) conditions relative to Schmidt and slow decrease (or oscillation) conditions relative to logarithmic summability in certain senses, respectively. Herein, the main issue to be discussed is in what ways the conditions imposed on the weights $(p_m), (q_n)$ or its partial sums $(P_m), (Q_n)$ should change while the ones imposed on the sequence $(u_{mn})$ become more inclusive. To reach an answer about this, we need the class $SV\hspace {-0.1 cm}A_{reg(\alpha)}$ and its characterization.

In the present paper, we are interested in relation between $(\overline{N},p,q)$ method, being product of relevant one-dimensional summability methods,  and $P$-convergence for double sequences. In accordance with this aim, we derive some Tauberian conditions, controlling $O_L$- and $O$-oscillatory behavior of a double sequence in certain senses, from $(\overline{N},p,q)$ summability to $P$- convergence with some restrictions on the weight sequences.

\section{Preliminaries}

In this section, we preface with basic definitions and notations in regards to double sequences and their weighted means.  Subsequent to these, we introduce the concepts of slow decrease relative to $(P_m)$ and $(Q_n),$ and slow oscillation relative to $(P_m)$ and $(Q_n)$ for double sequences of real and complex numbers and exhibit how a relation exists between  newly-described concepts. We put an end to this section by familiarizing the class $SV\hspace{-0.1 cm}A,$ its characterization and two of its subclasses.

A double sequence $u=(u_{mn})$ is a function $u$ from $\mathbb{N} \times \mathbb{N}$ into the set $\mathbb{K}$ ($\mathbb{K}$ is $\mathbb{R}$, the set of real numbers or $\mathbb{C}$, the set of complex numbers). The real or complex number $u_{mn}$ denotes the value of the function $u$ at a point $(m,n) \in \mathbb{N} \times \mathbb{N}$ and is called the $(m,n)$-term of the double sequence.\\
The set of all double sequences of real and complex numbers is denoted by $w^2(\mathbb{R})$ and $w^2(\mathbb{C})$, respectively.

A double sequence $(u_{mn})$ is said to be $P$-convergent to $\ell$ provided that for all $\epsilon >0$ there exists a $n_0=n_0(\epsilon)\in\mathbb{N}$ such that $|u_{mn}-\ell|<\epsilon$ whenever $m, n \geq n_0$ (see \cite{pringsheim}). The number $\ell$ is called the $P$-limit of $u$ and we denote it by $\displaystyle P-\hspace{-0.3cm}\lim _{m, n \to \infty} u_{mn}=\ell,$ where both $m$ and $n$ tend to $\infty$ independently of each other.\\
The set of all $P$-convergent double sequences of real and complex numbers is denoted by $c^2(\mathbb{R})$ and $c^2(\mathbb{C})$, respectively.

A double sequence $(u_{mn})$ is said to be bounded \textit{(or one-sided bounded)} provided that there exists a constant $M>0$ such that $|u_{mn}|\leq M$ \textit{(or $u_{mn}\geq -M$)} for all $m, n\in\mathbb{N}.$ \\
The set of all bounded double sequences of real and complex numbers is denoted by $\ell_{\infty}^{2}(\mathbb{R})$ and $\ell_{\infty}^{2}(\mathbb{C})$, respectively.\\
\indent Note that $(u_{mn})$ may converge without $(u_{mn})$ being a bounded function of $m$ and $n.$ To put it more explicitly, $P$-convergence of $(u_{mn})$ may not imply boundedness of its term in contrast to the case in single sequences. For example, the sequence $(u_{mn})$ defined by
$$u_{mn}=\begin{cases}
	7^{n} & \text{if} \;\; m=1; \; n=0,1,2,\dots\,, \\
	7^{m+2} & \text{if} \;\;n=3; \; m=0,1,2,\dots\,, \\
	2 & \text{otherwise}
\end{cases}$$
is $P$-convergent, but it is unbounded.

Some notations that will be used in places throughout this paper are given below.
\begin{notation}
	Let $(u_{mn})$ be a double sequence.
	\begin{itemize}
		\item[(i)] The symbol $u_{mn}=O(1)$ means that $|u_{mn}|\leq H $ for some constant $H>0$ and each $m, n\geq n_0.$
		\item[(ii)] The symbol $u_{mn}=O_L(1)$ means that $u_{mn}\geq M $ for some constant $M>0$ and each $m, n\geq n_0.$
		\item[(iii)] The symbol $u_{mn}=o(1)$ means that $u_{mn}\to 0$ as $m, n\to\infty.$
	\end{itemize}
\end{notation}

Let $u=(u_{mn})$ be a double sequence of real or complex numbers and let $p=(p_{mn})$ be a double sequence of positive integers such that
\begin{equation}\label{maincond}
	P_{mn}:=\sum_{i=0}^m\sum_{j=0}^{n}p_{ij}\to\infty\,\,\,\,\text{as}\,\,\,\,\max\{m,n\}\to\infty.
\end{equation}
The weighted means of $(u_{mn})$ with respect to the weights $(p_{mn})$ are defined by
\begin{equation}\label{weightedmeans_general}
	\sigma_{mn}:=\frac{1}{P_{mn}}\sum_{i=0}^m\sum_{j=0}^{n}p_{ij}u_{ij}
\end{equation}
for all $(m,n)\in\mathbb{N}\times\mathbb{N}$ and $P_{mn}>0.$

The following theorem proved by Kojima and Robinson \cite[Theorem 20]{hamilton} states necessary and sufficient conditions for regularity of transformation $\sigma_{mn}$ in the most general sense.
\begin{theorem}\label{KojimaRobinson}
	The necessary and sufficient conditions that every $P$-conver- gent double sequence imply $P$-convergence of its weighted means to same number under boundedness condition of double sequence are
	\begin{equation}\label{regularity-general}
		\frac{P_{mi}}{P_{mn}}\to 0\,\,\,\,\,\,\text{and}\,\,\,\,\,\,\frac{P_{jn}}{P_{mn}}\to 0\,\,\,\,\text{as}\,\,\,\,m, n\to\infty
	\end{equation}
	for any constant $i, j\in\mathbb{N}.$
\end{theorem}

In this paper, we deal only with a special class of weights which can be factorized in the form of $p_{mn}=p_{m}q_{n}$ where the single sequences $(p_{m})$ and $(q_{n})$ are positive weights such that
\begin{equation}\label{maincond-2}
	P_{m}:=\sum_{i=0}^{m} p_{i}\to\infty\,\,\,\,\,\,\,\,\text{and}\,\,\,\,\,\,\,\,Q_{n}:=\sum_{j=0}^{n} q_{j}\to\infty\,\,\,\,\text{as}\,\,\,\,m, n\to\infty.
\end{equation}
In the present case, the weighted means $(\sigma_{mn})$ defined in (\ref{weightedmeans_general}) transform $(u_{mn})$ into the form of
\begin{equation*}
	\sigma_{mn}:=\frac{1}{P_{m}Q_{n}}\sum_{i=0}^m\sum_{j=0}^{n}p_{i}q_{j}u_{ij}
\end{equation*}
for all $(m,n)\in\mathbb{N}\times\mathbb{N}$ and $P_{m}Q_{n}>0.$\\
A double sequence $(u_{mn})$ is called summable by the weighted mean method determined by the sequences $(p_{m})$ and $(q_{n})$ or shortly $(\overline{N},p,q)$ summable to $\ell$ provided that $(\sigma_{mn})$ is $P$-convergent to $\ell$. As a result of Theorem \ref{KojimaRobinson}, it can be seen that necessary and sufficient condition for regularity of the $(\overline{N},p,q)$ method is condition (\ref{maincond-2}). To put it another way, every $P$-convergent and bounded double sequence is also $(\overline{N},p,q)$ summable to same number under condition (\ref{maincond-2}). Nevertheless, the opposite of this proposition is not true in general. The question of whether some (nontrivial) condition on the terms $u_{mn}$ under which its $(\overline{N},p,q)$ summability implies its $P$-convergence exist comes to mind at this point. The condition $T\{u_{mn}\}$ making such a situation possible is called a \textit{Tauberian condition.} The resulting theorem stating that $P$-convergence follows from its $(\overline{N},p,q)$ summability and $T\{u_{mn}\}$ is called a \textit{Tauberian Theorem} for the $(\overline{N},p,q)$ method.

At present, we define the concepts of slow decrease and slow oscillation relative to $(P_m)$ and $(Q_n),$ for double sequences of real and complex numbers, respectively, besides we mention a relation between them.\\
A double sequence $(u_{mn})$ of real numbers is said to be slowly decreasing relative to both $(P_m)$ and $(Q_n)$ provided that
\begin{equation} \label{condbothSD}
	\lim_{\substack{\lambda\to 1^+\\\kappa\to 1^+}}\liminf_{m,n\to \infty}\min_{\substack{P_m\leq P_i\leq \lambda P_{m}\\
			Q_n\leq Q_j\leq \kappa Q_{n}}}\left(u_{ij}-u_{mn}\right)\geq0;
\end{equation}
that is, for each $\epsilon >0$ there exist $n_0=n_0(\epsilon)\in\mathbb{N},$ $\lambda =\lambda(\epsilon)>1$ and $\kappa =\kappa(\epsilon)>1$ such that
\begin{multline*}
	u_{ij}-u_{mn}\geq -\epsilon\,\,\,\text{whenever}\,\,\,n_0\leq m\leq i,\,\,\,n_0\leq n\leq j\,\,\,\text{and}\\
1\leq\frac{P_i}{P_m}\leq\lambda,\,\,\,1\leq\frac{Q_j}{Q_n}\leq\kappa.
\end{multline*}
Condition  (\ref{condbothSD}) is equivalent to
\begin{equation}
	\lim_{\substack{\lambda\to 1^-\\\kappa\to 1^-}}\liminf_{m,n\to \infty}\min_{\substack{ \lambda P_{m}<P_i\leq P_m\\
			\kappa Q_{n}<Q_j\leq Q_n}}\left(u_{mn}-u_{ij}\right)\geq0. \tag{\ref{condbothSD}$'$}
\end{equation}
A double sequence $(u_{mn})$ of complex numbers is said to be slowly oscillating relative to both $(P_m)$ and $(Q_n)$ provided that
\begin{equation} \label{condbothSO}
	\lim_{\substack{\lambda\to 1^+\\\kappa\to 1^+}}\limsup_{m,n\to \infty}\max_{\substack{P_m\leq P_i\leq \lambda P_{m}\\
			Q_n\leq Q_j\leq \kappa Q_{n}}}\left|u_{ij}-u_{mn}\right|=0;
\end{equation}
that is, for each $\epsilon >0$ there exist $n_0=n_0(\epsilon)\in\mathbb{N},$ $\lambda =\lambda(\epsilon)>1$ and $\kappa =\kappa(\epsilon)>1$ such that
\begin{multline*}
	|u_{ij}-u_{mn}|\leq\epsilon\,\,\,\text{whenever}\,\,\,n_0\leq m\leq i,\,\,\,n_0\leq n\leq j\,\,\,\text{and}\\
	1\leq\frac{P_i}{P_m}\leq\lambda,\,\,\,
	1\leq\frac{Q_j}{Q_n}\leq\kappa.
\end{multline*}
Condition  (\ref{condbothSO}) is equivalent to
\begin{equation*}
	\lim_{\substack{\lambda\to 1^-\\\kappa\to 1^-}}\limsup_{m,n\to \infty}\max_{\substack{ \lambda P_{m}<P_i\leq P_m\\
			\kappa Q_{n}<Q_j\leq Q_n}}\left|u_{mn}-u_{ij}\right|=0.\tag{\ref{condbothSO}$'$}
\end{equation*}
A double sequence $(u_{mn})$ of real numbers is said to be slowly decreasing relative to $(P_m)$ provided that
\begin{equation} \label{cond-1-SD}
	\lim_{\lambda\to 1^+}\liminf_{m,n\to \infty}\min_{P_m\leq P_i\leq \lambda P_{m}}\left(u_{in}-u_{mn}\right)\geq0,
\end{equation}
or equivalently,
\begin{equation}
	\lim_{\lambda\to 1^-}\liminf_{m,n\to \infty}\min_{\lambda P_{m}<P_i\leq P_m}\left(u_{mn}-u_{in}\right)\geq0,\tag{\ref{cond-1-SD}$'$}
\end{equation}
besides it is said to be slowly decreasing relative to $(P_m)$ in the strong sense if (\ref{cond-1-SD}) is satisfied with
\begin{equation} \label{cond-1-SD*}
	\min_{\substack{P_m\leq P_i\leq \lambda P_{m}\\
			Q_n\leq Q_j\leq \kappa Q_{n}}}\left(u_{ij}-u_{mj}\right)\,\,\,\,\,\text{instead of}\,\,\,\,\,\min_{P_m\leq P_i\leq \lambda P_{m}}\left(u_{in}-u_{mn}\right).
\end{equation}
A double sequence $(u_{mn})$ of complex numbers is said to be slowly oscillating relative to $(P_m)$ provided that
\begin{equation} \label{cond-1-SO}
	\lim_{\lambda\to 1^+}\limsup_{m,n\to \infty}\max_{P_m\leq P_i\leq \lambda P_{m}}\left|u_{in}-u_{mn}\right|=0,
\end{equation}
or equivalently,
\begin{equation}
	\lim_{\lambda\to 1^-}\limsup_{m,n\to \infty}\max_{\lambda P_{m}<P_i\leq P_m}\left|u_{mn}-u_{in}\right|=0,\tag{\ref{cond-1-SO}$'$}
\end{equation}
besides it is said to be slowly oscillating relative to $(P_m)$ in the strong sense if (\ref{cond-1-SO}) is satisfied with
\begin{equation} \label{cond-1-SO*}
	\max_{\substack{P_m\leq P_i\leq \lambda P_{m}\\
			Q_n\leq Q_j\leq \kappa Q_{n}}}\left|u_{ij}-u_{mj}\right|\,\,\,\,\,\text{instead of}\,\,\,\,\,\max_{P_m\leq P_i\leq \lambda P_{m}}\left|u_{in}-u_{mn}\right|.
\end{equation}
Similarly, the concepts of slow decrease and slow oscillation relative to $(Q_n)$ \textit{(in the strong sense)} for a double sequence $(u_{mn})$ of real and complex numbers  can be analogously defined, respectively.

Remark that if $(u_{mn})$ is slowly decreasing relative to $(P_m)$  in the strong sense and slowly decreasing relative to $(Q_n)$, then $(u_{mn})$ is slowly decreasing relative to both $(P_m)$ and $(Q_n).$\\
Indeed, for all large enough $m\,\,\text{and}\,\, n$, that is, $m,n\geq n_0,$ $\lambda>1,$ and $\kappa>1,$ we find
\begin{multline}\label{state1}
	\min_{\substack{P_m\leq P_i\leq \lambda P_{m}\\ Q_n\leq Q_j\leq \kappa Q_{n}}}\left(u_{ij}-u_{mn}\right)=
	\min_{\substack{P_m\leq P_i\leq \lambda P_{m}\\ Q_n\leq Q_j\leq \kappa Q_{n}}}\left(u_{ij}-u_{mj}+u_{mj}-u_{mn}\right)\nonumber\\
\geq\min_{\substack{P_m\leq P_i\leq \lambda P_{m}\\ Q_n\leq Q_j\leq \kappa Q_{n}}}\left(u_{ij}-u_{mj}\right)+
	\min_{\substack{ Q_n\leq j\leq \kappa Q_{n}}}\left(u_{mj}-u_{mn}\right).
\end{multline}
Taking lim inf and limit of both sides of (\ref{state1}) as $m,n\to\infty$ and $\lambda, \kappa \to 1^{+}$ respectively, we attain that the terms on the right-hand side of (\ref{state1}) are greater than $0$. Therefore, we reach that $(u_{mn})$ is slowly decreasing  relative to both $(P_m)$ and $(Q_n).$\\
In harmony with that, it can be said that if $(u_{mn})$ is slowly decreasing relative to $(P_m)$ and slowly decreasing relative to $(Q_n)$ in the strong sense, then it is slowly decreasing relative to both $(P_m)$ and $(Q_n).$

Similarly, if $(u_{mn})$ is slowly oscillating relative to $(P_m)$ and slowly oscillating relative to $(Q_n)$ in the strong sense, then $(u_{mn})$ is slowly oscillating relative to both $(P_m)$ and $(Q_n).$ In harmony with that, it can be said that if $(u_{mn})$ is slowly oscillating relative to $(Q_n)$ and slowly oscillating relative to $(P_m)$ in the strong sense, then it is slowly oscillating relative to both $(P_m)$ and $(Q_n).$

In the remainder of this section, we mention the classes including all positive sequences $(p_m)$ whose partial sum sequence $(P_m)$ is
\begin{itemize}
	\item[(i)] regularly varying sequence of positive index,
	\item[(ii)] rapidly varying sequence of index $\infty$  (see \cite{regularvariation} for more details).
\end{itemize}
Let $p=(p_m)$ be a sequence that satisfies $(p_m)=(P_m-P_{m-1}),$ where $P_{-1}=0$ and $P_m\neq 0$ for all $m\in\mathbb{N}.$
\begin{itemize}
	\item[(i)]  A sequence $(P_m)$ of positive numbers is said to be \textit{regularly varying} if for all $\lambda>0$
	\begin{equation}\label{regularlyvarying1}
		\lim_{m\to\infty}\frac{P_{\lambda_m}}{P_m}=\varphi(\lambda)\,\,\,\,\,\text{exists,}
	\end{equation}
	where $0<\varphi(\lambda)<\infty$ (cf. \cite{bojanicseneta}).\\
	\noindent In spite of the fact that this definition has been used by many authors as a starting point for studies including regularly varying sequences, these sequences possess quite useful properties, the most important of which is probably the following characterization theorem.
	\begin{namedtheorem}{Characterization Theorem} (\cite{karamata})
		The following statements are equivalent:
		\begin{itemize}
			\item[(a)] A sequence $(P_m)$ of positive numbers is a regularly varying sequence.
			\item[(b)] There exists a real number $\alpha>0$ such that $\varphi(\lambda)=\lambda^{\alpha}$ for all $\lambda>0.$
			\item[(c)] The sequence $(P_m)$ has the form $P_m=(m+1)^{\alpha}L(m)$ for $m\geq 0$ with constant $\alpha\geq 0$ and slowly varying function $L(.)$ on $(0, \infty),$  i.e. the function $L(.)$ is positive, measurable, and satisfies
			\begin{equation*}
				\lim_{t\to\infty}\frac{L(\lambda t)}{L(t)}=1\,\,\,\,\,\text{for all}\,\,\,\,\lambda>0.
			\end{equation*}
		\end{itemize}
	\end{namedtheorem}
	\noindent To emphasize such $\alpha,$ a sequence $(P_m)$ is called a \textit{regularly varying sequence of positive index $\alpha,$} as well. Note that a regularly varying sequence of index $\alpha=0$ corresponds to a slowly varying sequence.\\
	\noindent The set of all sequences of positive numbers $(p_m)$ satisfying $(c)$ is denoted by $SV\hspace {-0.1 cm}A_{reg(\alpha)}.$\\
	Here, it is useful to give the following implication proved by Bojanic and Seneta \cite{bojanicseneta}.
	\begin{lemma}\cite{bojanicseneta}\label{lemma2}
		If a sequence $P=(P_m)$ of positive numbers is regularly varying, then $P_{m-1}/P_m\to 1$ as $m\to\infty.$
	\end{lemma}
	\item[(ii)] A sequence $(P_m)$ of positive numbers is said to be \textit{rapidly varying of index $\infty$} if
	\begin{eqnarray}\label{rapidly}
		\frac{P_{\lambda_m}}{P_m}\to
		\left\{
		\begin{array}{lll}
			0 & \hbox{if $0<\lambda<1$,} \\
			1 & \hbox{if $\lambda=1$,} \\
			\infty & \hbox{if $\lambda>1$}
		\end{array}
		\right.\,\,\,\,\,\text{as}\,\,\,\,m\to\infty.
	\end{eqnarray}
	The set of all sequences of positive numbers $(p_m)$ satisfying (\ref{rapidly}) is denoted by $SV\hspace {-0.1 cm}A_{rap}.$\\
	In addition, it may be written conventionally as $\lambda^{\infty}$ because the right hand side of (\ref{rapidly}) is the limit of $\lambda^{\alpha}$ as $\alpha\to\infty.$
	%To establish (\ref{rapidly}), only $\displaystyle\frac{P_{\lambda_m}}{P_m}\to\infty$ for $\lambda>1$ has to be proved.
\end{itemize}

\section{Auxiliary results}\label{section3*}

In this section, we state an auxiliary result to be benefitted in the proofs of main results. The following lemma indicates two representations of difference between general terms of $(u_{mn})$ and $(\sigma_{mn})$ and it can be proved when it is make convenient modification in Lemma 1.2 which was presented by Fekete \cite{fekete}.
\begin{lemma}\label{lemma}
	Let $u=(u_{mn})$ be a double sequence.
	\begin{itemize}
		\item[(i)] For sufficiently large $\mu> m$ and $\eta>n,$ we have
		\begin{eqnarray*}
			u_{mn}-\sigma_{mn}&=&\frac{P_{\mu}Q_{\eta}}{(P_{\mu}-P_m)(Q_{\eta}-Q_n)}\left(\sigma_{\mu\eta}-\sigma_{\mu n}-\sigma_{m\eta}+\sigma_{mn}\right)\\
			&+&\frac{P_{\mu}}{P_{\mu}-P_m}\left(\sigma_{\mu n}-\sigma_{mn}\right)
			+\frac{Q_{\eta}}{Q_{\eta}-Q_n}\left(\sigma_{m\eta}-\sigma_{mn}\right)\\
			&-&\frac{1}{(P_{\mu}-P_m)(Q_{\eta}-Q_n)}		 \sum_{i=m+1}^{\mu}\sum_{j=n+1}^{\eta}p_iq_j\left(u_{ij}-u_{mn}\right).
		\end{eqnarray*}
		\item[(ii)] For sufficiently large $\mu< m$ and $\eta<n,$ we have
		\begin{eqnarray*}
			u_{mn}-\sigma_{mn}&=&\frac{P_{\mu}Q_{\eta}}{(P_m-P_{\mu})(Q_n-Q_{\eta})}\left(\sigma_{mn}-\sigma_{\mu n}-\sigma_{m\eta}+\sigma_{\mu\eta}\right)\\
			&+&\frac{P_{\mu}}{P_m-P_{\mu}}\left(\sigma_{mn}-\sigma_{\mu n}\right)+\frac{Q_{\eta}}	{Q_n-Q_{\eta}}\left(\sigma_{mn}-\sigma_{m\eta}\right)\\
			&+&\frac{1}{(P_m-P_{\mu})(Q_n-Q_{\eta})}	 \sum_{i=\mu+1}^{m}\sum_{j=\eta+1}^{n}p_iq_j\left(u_{mn}-u_{ij}\right).
		\end{eqnarray*}
	\end{itemize}
\end{lemma}

\section{Main Results for the $(\overline{N}, p, q)$ Summable Double Sequences of real numbers}\label{section3}

%\subsection{Main Results}\label{section3.1}\hspace*{\fill} \\

%\subsection{Examples}\label{section3.2}\hspace*{\fill} \\

%\subsection{Proofs of Main Results}\label{section3.3}\hspace*{\fill} \\

This section is constructed by considering the following headings for $(\overline{N},p,q)$ summable double sequences of real numbers:
\begin{itemize}
	\item[(a)] Determining certain subsets $\mathcal{T}\{u_{mn}\}$ of double sequence space and certain subclasses $\mathcal{C}\{p_m, q_n\}$ of single sequence space having the property: ``Under conditions  $\mathcal{C}\{p_m, q_n\},$ a double sequence $(u_{mn})\in \mathcal{T}\{u_{mn}\}$ being $(\overline{N},p,q)$ summable to $l\in \mathbb{R}$ is also $P$-convergent to same value.''
	\item[(b)] Demonstrating with examples how the subsets $\mathcal{T}\{u_{mn}\}$ and the subclasses $\mathcal{C}\{p_m, q_n\}$ change for special means occurring depends on choosing of weight sequences $(p_{m})$ and $(q_{n}).$
	\item[(c)] Proving the property given in (a) for the related subsets and subclasses.
\end{itemize}
Here, we formulate our main results for double sequences of real numbers as follows:
\begin{theorem}\label{maintheo-1}
	Let $(p_m), (q_n)\in SV\hspace {-0.1 cm}A_{reg(\alpha)}$ and a double sequence $(u_{mn})$ be $(\overline{N},p,q)$ summable to a number $\ell.$
	If $(u_{mn})$ is slowly decreasing relative to $(P_m),$ slowly decreasing relative to $(Q_n),$ and slowly decreasing relative to $(P_m)$ or $(Q_n)$ in the strong sense, then $(u_{mn})$ is $P$-convergent to $\ell.$
\end{theorem}

\begin{theorem}\label{maintheo-2}
	Let $(p_m), (q_n)\in SV\hspace {-0.1 cm}A_{reg(\alpha)}$ and a double sequence $(u_{mn})$ be $(\overline{N},p,q)$ summable to a number $\ell.$  If $(u_{mn})$ satisfies conditions
	\begin{equation}\label{landautype(P)}
		\frac{P_{m}}{p_m}\Delta_{10}u_{mn}=O_L(1)\,\,\,\,\,\text{and}\,\,\,\,\,\frac{Q_{n}}{q_n}\Delta_{01}u_{mn}=O_L(1),
	\end{equation}
	then $(u_{mn})$ is $P$-convergent to $\ell.$	
\end{theorem}

\begin{exp-non}
	In conjunction with the weighted means, there are many special means occurring depends on choosing of weight sequences $(p_{m})$ and $(q_{m})$.  Included by the weighted means and also commonly used by researchers in literature, some means are listed with their corresponding subsets $\mathcal{T}\{u_{mn}\}$ and subclasses $\mathcal{C}\{p_m, q_n\}$ for double sequences of real numbers as follows.
	\begin{itemize}
		\item[(i)]  In case $p_{m}=q_n=1,$ it leads to the \textit{arithmetic means} (or called the Ces\`{a}ro means of order $(1,1)$) of a double sequence where $P_{m}Q_n=(m+1)(n+1)$ for all $m, n\in\mathbb{N}.$ Under the circumstances, the conditions in Theorem \ref{maintheo-1} correspond to slow decrease conditions in senses $(1,0), (0,1)$ and in the strong sense $(1,0)$ or $(0,1),$ besides conditions (\ref{landautype(P)}) correspond to conditions
		\begin{equation*}
			m\Delta_{10}u_{mn}=O_L(1)\,\,\,\,\,\text{and}\,\,\,\,\,n\Delta_{01}u_{mn}=O_L(1),
		\end{equation*}
		however slow decrease condition in sense $(1,1)$ and condition $mn\Delta_{11}$ $u_{mn}=O_L(1)$ discussed by M\'{o}ricz  \cite{moricz1994} therein are superfluous. As a result of our main Theorem \ref{maintheo-1} and Theorem \ref{maintheo-2}, the mentioned conditions are sufficient Tauberian conditions for the $(C,1,1)$ summability of double sequences of real numbers.
		\item[(ii)] In case $\displaystyle{p_{m}q_n=\frac{1}{(m+1)(n+1)}},$ it leads to the \textit{harmonic means} (or called the logarithmic means) of a double sequence where $\displaystyle{P_{m}Q_n\sim}$ $\log m\log n$ for all $m, n\in\mathbb{N}.$ Under the circumstances, the conditions in Theorem \ref{maintheo-1} correspond to conditions of slow decrease with respect to summability $(L,1)$ in senses $(1,0), (0,1)$ in the strong sense $(1,0)$ or $(0,1),$i.e.,
		\begin{equation}\label{SDlog10}
			\lim_{\lambda\to 1^+}\liminf_{m,n\to \infty}\min_{m\leq i\leq m^\lambda}(u_{in}-u_{mn})\geq 0,
		\end{equation}
		\begin{equation}\label{SDlog01}
			\lim_{\kappa\to 1^+}\liminf_{m,n\to \infty}\min_{n\leq j\leq n^\kappa}(u_{mj}-u_{mn})\geq 0
		\end{equation}
		and shall we say
		\begin{equation}\label{SDlog01strong}
			\lim_{\lambda, \kappa\to 1^+}\liminf_{m,n\to \infty}\min_{\substack{m\leq i\leq m^\lambda\\ n\leq j\leq n^\kappa}}(u_{ij}-u_{in})\geq 0,
		\end{equation}
		which was discussed by Kwee \cite{kwee} and M\'{o}ricz \cite{moricz2013} for the logarithmic summability of single sequences in different ways. In addition, conditions (\ref{landautype(P)}) correspond to conditions
		\begin{equation}\label{landautypelog}
			m\log (m+1)\Delta_{10}u_{mn}=O_L(1)\,\,\,\,\,\text{and}\,\,\,\,\,n\log (n+1)\Delta_{01}u_{mn}=O_L(1).
		\end{equation}
		As a result of our main Theorem \ref{maintheo-1} and Theorem \ref{maintheo-2}, the mentioned conditions are sufficient Tauberian conditions for the logarithmic summability of double sequences of real numbers.
	\end{itemize}
\end{exp-non}

Now, we present the proofs of our main results for double sequences of real numbers.

\noindent \textbf{\emph{Proof of Theorem} \ref{maintheo-1}}.
Assume that $(u_{mn})$ being $(\overline{N},p,q)$ summable to $\ell$ is slowly decreasing relative to $(P_m)$ and $(Q_n)$ and slowly decreasing relative to $(Q_n)$ in the strong sense. In order to prove that $(u_{mn})$ is $P$-convergent to same number, we indicate that difference $u_{mn}-\sigma_{mn}$ is $P$-convergent to $0.$ By means of Lemma \ref{lemma}, we have for $\mu> m$ and $\eta>n$
\begin{eqnarray}\label{maintheo-eq1}
	u_{mn}-\sigma_{mn}&=&\frac{P_{\mu}Q_{\eta}}{(P_{\mu}-P_m)(Q_{\eta}-Q_n)}\left(\sigma_{\mu\eta}-\sigma_{\mu n}-\sigma_{m\eta}+\sigma_{mn}\right)\nonumber\\
	&+&\frac{P_{\mu}}{P_{\mu}-P_m}\left(\sigma_{\mu n}-\sigma_{mn}\right)
	+\frac{Q_{\eta}}{Q_{\eta}-Q_n}\left(\sigma_{m\eta}-\sigma_{mn}\right)\nonumber\\
	&-&\frac{1}{(P_{\mu}-P_m)(Q_{\eta}-Q_n)}\sum_{i=m+1}^{\mu}\sum_{j=n+1}^{\eta}p_iq_j\left(u_{ij}-u_{mn}\right)\nonumber\\
	&\leq& \frac{P_{\mu}Q_{\eta}}{(P_{\mu}-P_m)(Q_{\eta}-Q_n)}\left(\sigma_{\mu\eta}-\sigma_{\mu n}-\sigma_{m\eta}+\sigma_{mn}\right)\nonumber\\
	&+&\frac{P_{\mu}}{P_{\mu}-P_m}\left(\sigma_{\mu n}-\sigma_{mn}\right)
	+\frac{Q_{\eta}}{Q_{\eta}-Q_n}\left(\sigma_{m\eta}-\sigma_{mn}\right)\nonumber\\
	&-&\min_{\substack{m\leq i\leq \mu\\ n\leq j\leq \eta}}	\left(u_{ij}-u_{in}\right)-\min_{m\leq i\leq\mu}\left(u_{in}-u_{mn}\right).
\end{eqnarray}
Putting $$\mu=\argmin\left\{P_i\geq\left(1+\frac{\delta}{2}\right) P_m \right\}=\min\left\{i> m:P_i\geq\left(1+\frac{\delta}{2}\right) P_m\right\}$$ and $$\eta=\argmin\left\{Q_j\geq\left(1+\frac{\gamma}{2}\right)  Q_n\right\}=\min\left\{j> n: Q_j\geq\left(1+\frac{\gamma}{2}\right) Q_n\right\}$$ with $\delta, \gamma>0,$ we get $\mu\geq m,\,\eta\geq n$ and $\displaystyle P_{\mu-1}<\left(1+\frac{\delta}{2}\right) P_m$ and $\displaystyle Q_{\eta-1}<\left(1+\frac{\gamma}{2}\right) Q_n.$ Moreover, on account of $(p_m), (q_n)\in SV\hspace {-0.1 cm}A_{reg(\alpha)},$ we have $\displaystyle \frac{P_m}{P_{m-1}}$ $\displaystyle \leq 1+\frac{\delta}{4}$ and $\displaystyle \frac{Q_n}{Q_{n-1}}\leq 1+\frac{\gamma}{4}$ for sufficiently large $m, n.$ By means of Lemma \ref{lemma2}, it can be easily seen that $p_m/P_m\to 0$ and $q_n/Q_n\to 0$ as $m, n\to\infty,$ as well. As a result of these findings, we obtain by Lemma \ref{lemma2} that
\begin{equation}\label{O(Pm)}
	P_\mu=\frac{P_\mu}{P_{\mu-1}}P_{\mu-1}\leq\left(1+\frac{\delta}{4}\right)P_m\left(1+\frac{\delta}{2}\right)\leq P_m(1+\delta)=\lambda P_m,
\end{equation}
or in another saying,
\begin{equation}\label{O(Pm)*}
	\frac{P_\mu}{P_m}=\frac{P_{\mu-1}}{P_m}+\frac{p_\mu}{P_m}\leq\lambda+\frac{p_\mu}{P_\mu}\frac{P_\mu}{P_m}=\lambda(1+o(1)),\,\,\,\,\,\,\,\,
	\frac{P_\mu}{P_m}\geq\left(1+\frac{\delta}{2}\right)=\frac{\lambda+1}{2}
\end{equation}
and
\begin{equation}\label{O(Qn)}
	Q_\eta=\frac{Q_\eta}{Q_{\eta-1}}Q_{\eta-1}\leq\left(1+\frac{\gamma}{4}\right)Q_n\left(1+\frac{\gamma}{2}\right)\leq Q_n(1+\gamma)=\kappa Q_n
\end{equation}
or in another saying,
\begin{equation}\label{O(Qn)*}
	\frac{Q_\eta}{Q_n}=\frac{Q_{\eta-1}}{Q_n}+\frac{q_\eta}{Q_n}\leq\kappa+\frac{q_\eta}{Q_\eta}\frac{Q_\eta}{Q_n}=\kappa(1+o(1))\,\,\,\,\,\,\,\,
	\frac{Q_\eta}{Q_n}\geq\left(1+\frac{\gamma}{2}\right)=\frac{\kappa+1}{2}
\end{equation}
for $\lambda,\,\kappa>1.$
%Putting
%$$\mu=\argmin\{P_i\geq\lambda P_m \}=\min\{i> m:P_i\geq\lambda P_m\}$$ and $$\eta=\argmin\{Q_j\geq\kappa Q_n\}=\min\{j> n: Q_j\geq\kappa Q_n\}$$ with %$\lambda, \kappa>1,$  we observe by Lemma \ref{lemma2} that $p_m/P_m\to 0$ and $q_n/Q_n\to 0$ as $m, n\to\infty$ and hence
%\begin{equation}\label{O(Pm)}
	%\frac{P_\mu}{P_m}\geq\frac{\lambda P_m}{P_m}=\lambda\,\,\,\,\,\,\,\,\text{and}\,\,\,\,\,\,\,\,\frac{P_\mu}{P_m}=\frac{P_{\mu-1}}{P_m}+\frac{p_\mu}{P_m}
	%\leq\lambda+\frac{p_\mu}{P_\mu}\frac{P_\mu}{P_m}=\lambda(1+o(1))
	%\end{equation}
%and
%\begin{equation}\label{O(Qn)}
	%\frac{Q_\eta}{Q_n}\geq\frac{\kappa Q_n}{Q_n}=\kappa\,\,\,\,\,\,\,\,\text{and}\,\,\,\,\,\,\,\,\frac{Q_\eta}{Q_n}=\frac{Q_{\eta-1}}{Q_n}+\frac{q_\eta}{Q_n}
	%\leq\kappa+\frac{q_\eta}{Q_\eta}\frac{Q_\eta}{Q_n}=\kappa(1+o(1)),
	%\end{equation}
%which mean
%\begin{equation}\label{O(Pm)O(Qn)}
	%\frac{P_\mu}{P_m}\to\lambda\,\,\,\,\,\,\,\,\text{and}\,\,\,\,\,\,\,\,\frac{Q_\eta}{Q_n}\to\kappa\,\,\,\text{as}\,\,\, m,n\to\infty,
	%\end{equation}
%respectively.
From the assumption, we have $\sigma_{mn}\to\ell$ and so, $\sigma_{ij}-\sigma_{mn}\to 0$ for $i=m$ or $\mu$ and $j=n$ or $\eta$ as $m, n\to\infty.$
Since the sequences $(P_m)$ and $(Q_n)$ are strictly increasing sequences, we attain from (\ref{maintheo-eq1}) that
\begin{eqnarray}\label{maintheo-eq2}
	u_{mn}-\sigma_{mn} &\leq& \frac{P_{\mu}Q_{\eta}}{(P_{\mu}-P_m)(Q_{\eta}-Q_n)}\left(\sigma_{\mu\eta}-\sigma_{\mu n}-\sigma_{m\eta}+\sigma_{mn}\right)\nonumber\\
	&+&\frac{P_{\mu}}{P_{\mu}-P_m}\left(\sigma_{\mu n}-\sigma_{mn}\right)
	+\frac{Q_{\eta}}{Q_{\eta}-Q_n}\left(\sigma_{m\eta}-\sigma_{mn}\right)\nonumber\\
	&-&\min_{\substack{P_m\leq P_i\leq \lambda P_{m}\\ Q_n\leq Q_j\leq \kappa Q_{n}}}
	\left(u_{ij}-u_{in}\right)-\min_{P_m\leq P_i\leq \lambda P_{m}}\left(u_{in}-u_{mn}\right)\nonumber\\
	&\leq& \frac{4\lambda\kappa+2\lambda(\kappa-1)+2(\lambda-1)\kappa}{(\lambda-1)(\kappa-1)}(1+o(1)) o(1)\nonumber\\
	&-&\min_{\substack{P_m\leq P_i\leq \lambda P_{m}\\ Q_n\leq Q_j\leq \kappa Q_{n}}}
	\left(u_{ij}-u_{in}\right)-\min_{P_m\leq P_i\leq \lambda P_{m}}\left(u_{in}-u_{mn}\right)
\end{eqnarray}
where
\begin{multline*}
	\frac{P_\mu}{P_\mu-P_m}=\frac{P_\mu/P_m}{P_\mu/P_m-1}\leq\frac{2\lambda}{(\lambda-1)}(1+o(1))\,\,\,\,\,\,\,\,\text{and}\\
	\frac{Q_\eta}{Q_\eta-Q_n}=\frac{Q_\eta/Q_n}{Q_\eta/Q_n-1}\leq\frac{2\kappa}{(\kappa-1)}(1+o(1)).
\end{multline*}
If we get $\limsup$ of both sides of inequality (\ref{maintheo-eq2}) as $m, n\to\infty,$ then we reach for any $\lambda, \kappa>1$
\begin{multline*}
	\limsup_{m, n\to\infty} (u_{mn}-\sigma_{mn})\leq -\liminf_{m, n\to\infty}\min_{\substack{P_m\leq P_i\leq \lambda P_{m}\\ Q_n\leq Q_j\leq \kappa Q_{n}}}
	\left(u_{ij}-u_{in}\right)\\
	-\liminf_{m, n\to\infty}\min_{P_m\leq P_i\leq \lambda P_{m}}\left(u_{in}-u_{mn}\right).
\end{multline*}
If we get limit of both sides of last inequality as $\lambda, \kappa\to 1^+,$ then we find
\begin{equation}\label{maintheo-eq3}
	\limsup_{m, n\to\infty} (u_{mn}-\sigma_{mn})\leq 0
\end{equation}
due to the fact that $(u_{mn})$ is slowly decreasing relative to $(P_m)$ and $(Q_n)$ and slowly decreasing relative to $(Q_n)$ in the strong sense. Following a similar procedure to above for $\tilde{\mu}< m$ and $\tilde{\eta}<n,$ we have by means of Lemma \ref{lemma}
\begin{eqnarray}\label{maintheo-eq11}
	u_{mn}-\sigma_{mn}&=&\frac{P_{\tilde{\mu}}Q_{\tilde{\eta}}}{(P_m-P_{\tilde{\mu}})(Q_n-Q_{\tilde{\eta}})}\left(\sigma_{mn}-\sigma_{\tilde{\mu} n}-\sigma_{m\tilde{\eta}}+\sigma_{\tilde{\mu}\tilde{\eta}}\right)\nonumber\\
	&+&\frac{P_{\tilde{\mu}}}{P_m-P_{\tilde{\mu}}}\left(\sigma_{mn}-\sigma_{\tilde{\mu} n}\right)
	+\frac{Q_{\tilde{\eta}}}	{Q_n-Q_{\tilde{\eta}}}\left(\sigma_{mn}-\sigma_{m\tilde{\eta}}\right)\nonumber\\
	&+&\frac{1}{(P_m-P_{\tilde{\mu}})(Q_n-Q_{\tilde{\eta}})}	 \sum_{i=\tilde{\mu}+1}^{m}\sum_{j=\tilde{\eta}+1}^{n}p_iq_j\left(u_{mn}-u_{ij}\right)\nonumber\\
	&\geq& \frac{P_{\tilde{\mu}}Q_{\tilde{\eta}}}{(P_m-P_{\tilde{\mu}})(Q_n-Q_{\tilde{\eta}})}\left(\sigma_{mn}-\sigma_{\tilde{\mu} n}-\sigma_{m\tilde{\eta}}+\sigma_{\tilde{\mu}\tilde{\eta}}\right)\nonumber\\
	&+&\frac{P_{\tilde{\mu}}}{P_m-P_{\tilde{\mu}}}\left(\sigma_{mn}-\sigma_{\tilde{\mu} n}\right)
	+\frac{Q_{\tilde{\eta}}}	{Q_n-Q_{\tilde{\eta}}}\left(\sigma_{mn}-\sigma_{m\tilde{\eta}}\right)\nonumber\\
	&+&\min_{\tilde{\mu}\leq i\leq m}\left(u_{mn}-u_{in}\right)+\min_{\substack{\tilde{\mu}\leq i\leq m\\ \tilde{\eta}\leq j\leq n}}\left(u_{in}-u_{ij}\right).
\end{eqnarray}
Putting $$\tilde{\mu}=\argmax\left\{P_m\geq\left(1+\frac{\delta}{2}\right)P_i \right\}=\max\left\{m>i: P_m\geq\left(1+\frac{\delta}{2}\right)P_i \right\}$$ and $$\tilde{\eta}=\argmax\left\{Q_n\geq\left(1+\frac{\gamma}{2}\right)Q_j\right\}=\max\left\{n>j : Q_n\geq\left(1+\frac{\delta}{2}\right)Q_j\right\}$$  with $\delta, \gamma>0,$  we get $\tilde{\mu}\leq m,\,\tilde{\eta}\leq n$ and $\displaystyle \left(1+\frac{\delta}{2}\right)P_{\tilde{\mu}+1}>P_m$ and $\displaystyle \left(1+\frac{\gamma}{2}\right)$ $Q_{\tilde{\eta}+1}>Q_n.$ Moreover, on account of $(p_m), (q_n)\in SV\hspace {-0.1 cm}A_{reg(\alpha)},$ we have
$\displaystyle \frac{P_{\tilde{\mu}+1}}{P_{\tilde{\mu}}}\leq 1+\frac{\delta}{4}$ and $\displaystyle \frac{Q_{\tilde{\eta}+1}}{Q_{\tilde{\eta}}}\leq 1+\frac{\gamma}{4}$ for sufficiently large $\tilde{\mu}, \tilde{\eta}.$ As a result of these findings, we obtain by Lemma \ref{lemma2} that
\begin{multline*}
	P_{\tilde{\mu}}=\frac{P_{\tilde{\mu}}}{P_{\tilde{\mu}+1}}P_{\tilde{\mu}+1}\geq \frac{1}{\displaystyle\left(1+\frac{\delta}{4}\right)}\frac{1}{\displaystyle\left(1+\frac{\delta}{2}\right)}P_m\geq\frac{1}{\displaystyle\left(1+\frac{\delta}{2}\right)^2}P_m=\tilde{\lambda}P_m\,\,\,\,\,\,\,\,\text{and}\\
\frac{P_{\tilde{\mu}}}{P_m}\leq\frac{1}{\displaystyle\left(1+\frac{\delta}{2}\right)}=\sqrt{\tilde{\lambda}}
\end{multline*}
and
\begin{multline*}
	Q_{\tilde{\eta}}=\frac{Q_{\tilde{\eta}}}{Q_{\tilde{\eta}+1}}Q_{\tilde{\eta}+1}\geq \frac{1}{\displaystyle\left(1+\frac{\gamma}{4}\right)}\frac{1}{\displaystyle\left(1+\frac{\gamma}{2}\right)}Q_n\geq\frac{1}{\displaystyle\left(1+\frac{\gamma}{2}\right)^2}Q_n=\tilde{\kappa}Q_n\,\,\,\,\,\,\,\,\text{and}\\
\frac{Q_{\tilde{\eta}}}{Q_n}\leq\frac{1}{\displaystyle\left(1+\frac{\gamma}{2}\right)}=\sqrt{\tilde{\kappa}}
\end{multline*}
for $0<\tilde{\lambda},\,\tilde{\kappa}<1.$
%Putting $$\tilde{\mu}=\argmax\{P_m\geq\lambda P_i \}=\max\{m>i: P_m\geq\lambda P_i \}$$ and $$\tilde{\eta}=\argmax\{Q_n\geq\kappa Q_j\}=\max\{n>j : %Q_n\geq\kappa Q_j\}$$ with $\lambda, \kappa>1,$  we observe by Lemma \ref{lemma2} that $p_{\tilde{\mu}}/P_{\tilde{\mu}}\to 0$ and %$q_{\tilde{\eta}}/Q_{\tilde{\eta}}\to 0$ as $\tilde{\mu}, \tilde{\eta}\to\infty$ and hence
%\begin{equation*}
	%\frac{P_{\tilde{\mu}}}{P_m}\leq\frac{1}{\lambda}\,\,\,\,\,\,\,\,\text{and}\,\,\,\,\,\,\,\,
	%\frac{P_{\tilde{\mu}}}{P_m}=\frac{P_{\tilde{\mu}+1}}{P_m}-\frac{p_{\tilde{\mu}+1}}{P_m}>\frac{1}{\lambda}-
	%\frac{p_{\tilde{\mu}+1}}{P_{\tilde{\mu}+1}}\frac{P_{\tilde{\mu}+1}}{P_{\tilde{\mu}}}\frac{{P_{\tilde{\mu}}}}{P_m}=\frac{1}{\lambda}-\frac{1}{\lambda}o(1)(1+o(1))=\frac{1}{%\lambda}(1-o(1))
	%\end{equation*}
%and
%\begin{equation*}
	%	\frac{Q_{\tilde{\eta}}}{Q_n}\leq\frac{1}{\kappa}\,\,\,\,\,\,\,\,\text{and}\,\,\,\,\,\,\,\,
	%	\frac{Q_{\tilde{\eta}}}{Q_n}=\frac{Q_{\tilde{\eta}+1}}{Q_n}-\frac{q_{\tilde{\eta}+1}}{Q_n}>\frac{1}{\kappa}-
	%	\frac{q_{\tilde{\eta}+1}}{Q_{\tilde{\eta}+1}}\frac{Q_{\tilde{\eta}+1}}{Q_{\tilde{\eta}}}\frac{{Q_{\tilde{\eta}}}}{Q_n}=\frac{1}{\kappa}-\frac{1}{\kappa}
	%	o(1)(1+o(1))=\frac{1}{\kappa}(1-o(1))
	%\end{equation*}
%which mean
%\begin{equation*}
	%\frac{P_{\tilde{\mu}}}{P_m}\to\frac{1}{\lambda}\,\,\,\,\,\,\,\,\text{and}\,\,\,\,\,\,\,\,\frac{Q_{\tilde{\eta}}}{Q_n}\to\frac{1}{\kappa} \,\,\,\text{as}\,\,\, \tilde{\mu}, %\tilde{\eta}\to\infty,
	%\end{equation*}
%respectively.
From the assumption, we have $\sigma_{mn}\to\ell$ and so, $\sigma_{mn}-\sigma_{ij}\to 0$ for $i=m$ or $\tilde{\mu}$ and $j=n$ or $\tilde{\eta}$ as $\tilde{\mu}, \tilde{\eta}\to\infty.$ Since the sequences $(P_m)$ and $(Q_n)$ are strictly increasing sequences, we attain from (\ref{maintheo-eq11}) that
\begin{eqnarray}\label{maintheo-eq21}
	u_{mn}-\sigma_{mn}&\geq& \frac{P_{\tilde{\mu}}Q_{\tilde{\eta}}}{(P_m-P_{\tilde{\mu}})(Q_n-Q_{\tilde{\eta}})}\left(\sigma_{mn}-\sigma_{\tilde{\mu} n}-\sigma_{m\tilde{\eta}}+\sigma_{\tilde{\mu}\tilde{\eta}}\right)\nonumber\\
	&+&\frac{P_{\tilde{\mu}}}{P_m-P_{\tilde{\mu}}}\left(\sigma_{mn}-\sigma_{\tilde{\mu} n}\right)
	+\frac{Q_{\tilde{\eta}}}	{Q_n-Q_{\tilde{\eta}}}\left(\sigma_{mn}-\sigma_{m\tilde{\eta}}\right)\nonumber\\
	&+&\min_{\tilde{\lambda} P_m<P_i\leq P_m}\left(u_{mn}-u_{in}\right)+\min_{\substack{\tilde{\lambda} P_{m}\leq P_i\leq P_{m}\\ \tilde{\kappa} Q_{n}\leq Q_j\leq  Q_{n}}}\left(u_{in}-u_{ij}\right)\nonumber\\
	&\geq&\frac{\tilde{\lambda}+\tilde{\kappa}-\tilde{\lambda}\tilde{\kappa}}{(1-\tilde{\lambda})(1-\tilde{\kappa})}o(1)+\min_{\tilde{\lambda} P_m<P_i\leq P_m}\left(u_{mn}-u_{in}\right)\nonumber\\
	&+&\min_{\substack{\tilde{\lambda} P_{m}\leq P_i\leq P_{m}\\ \tilde{\kappa} Q_{n}\leq Q_j\leq  Q_{n}}}\left(u_{in}-u_{ij}\right)
\end{eqnarray}
where
\begin{equation*}
	 \frac{P_{\tilde{\mu}}}{P_m-P_{\tilde{\mu}}}=\frac{P_{\tilde{\mu}}/P_m}{1-P_{\tilde{\mu}}/P_m}\geq\frac{\tilde{\lambda}}{1-\tilde{\lambda}}\,\,\,\,\,\text{and}\,\,\,\,\,\frac{Q_{\tilde{\eta}}}{Q_n-Q_{\tilde{\eta}}}=\frac{Q_{\tilde{\eta}}/Q_n}{1-Q_{\tilde{\eta}}/Q_n}\geq\frac{\tilde{\kappa}}{1-\tilde{\kappa}},
\end{equation*}	
for $1/\lambda=\tilde{\lambda},$ $1/\kappa=\tilde{\kappa},$ and $0<\tilde{\lambda}, \tilde{\kappa}<1.$ If we get $\liminf$ of both sides of inequality (\ref{maintheo-eq21}) as $\tilde{\mu}, \tilde{\eta}\to\infty,$ then we reach for any $0<\tilde{\lambda}, \tilde{\kappa}<1$
\begin{multline*}
	\liminf_{m, n\to\infty} (u_{mn}-\sigma_{mn})\geq\liminf_{m, n\to\infty}\min_{\tilde{\lambda} P_m<P_i\leq P_m}\left(u_{mn}-u_{in}\right)\\
	+\liminf_{m, n\to\infty}\min_{\substack{\tilde{\lambda} P_{m}\leq P_i\leq P_{m}\\ \tilde{\kappa} Q_{n}\leq Q_j\leq  Q_{n}}}\left(u_{in}-u_{ij}\right).
\end{multline*}
If we get limit of both sides of last inequality as $\tilde{\lambda}, \tilde{\kappa}\to 1^-,$ then we find
\begin{equation}\label{maintheo-eq31}
	\liminf_{m, n\to\infty} (u_{mn}-\sigma_{mn})\geq 0
\end{equation}
due to the fact that $(u_{mn})$ is slowly decreasing relative to $(P_m)$ and $(Q_n)$ and slowly decreasing relative to $(Q_n)$ in the strong sense. If we combine inequalities (\ref{maintheo-eq3}) with (\ref{maintheo-eq31}), we conclude
\begin{equation*}
	\lim_{m, n\to\infty} (u_{mn}-\sigma_{mn})= 0,
\end{equation*}
which means that $(u_{mn})$ is $P$-convergent to $\ell.$ \qed\\

\noindent \textbf{\emph{Proof of Theorem} \ref{maintheo-2}}.
Assume that $(u_{mn})$ being $(\overline{N},p,q)$ summable to $\ell$ satisfies conditions (\ref{landautype(P)}). If we indicate that conditions (\ref{landautype(P)}) imply slow decrease condition relative to $(P_m)$ and $(Q_n)$ and slow decrease condition relative to $(Q_n)$ in the strong sense, then we could this proof with the help of Theorem \ref{maintheo-1}. \\
Putting $\mu=\argmin\left\{P_i\geq\left(1+\displaystyle\frac{\delta}{2}\right) P_m \right\}$ and $\eta=\argmin\left\{Q_j\geq\left(1+\displaystyle\frac{\gamma}{2}\right)  Q_n\right\}$ with $\delta, \gamma>0,$ we can observe inequalities  (\ref{O(Pm)*})  and (\ref{O(Qn)*}). Then, we have for $n_0\leq m\leq i\leq\mu$ and $n_0\leq n$
\begin{eqnarray*}
	u_{in}-u_{mn}=\sum_{k=m+1}^{i}\Delta_{10}u_{kn}\geq -M_1\sum_{k=m+1}^{i}\frac{p_k}{P_k}&\geq& -M_1\left(\frac{P_\mu}{P_m}-1\right)\\
	&\geq& -M_1(\lambda-1+\lambda o(1))
\end{eqnarray*}
for any constant $M_1>0$ and $\lambda>1.$ If we get $\liminf$ and limit of both sides of last inequality as $m, n\to\infty$ and $\lambda\to 1^+$ respectively, then we reach
\begin{equation*}
	\lim_{\lambda\to 1^+}\liminf_{m,n\to \infty}\min_{P_m\leq P_i\leq \lambda P_{m}}\left(u_{in}-u_{mn}\right)\geq0,
\end{equation*}
which means that $(u_{mn})$ is slowly decreasing relative to $(P_m)$. Similarly, we obtain $n_0\leq n\leq j\leq\eta$ and $n_0\leq m$
\begin{eqnarray*}
	u_{mj}-u_{mn}=\sum_{r=n+1}^{j}\Delta_{01}u_{mr}\geq -M_2\sum_{r=n+1}^{j}\frac{q_r}{Q_r}&\geq& -M_2\left(\frac{Q_\eta}{Q_n}-1\right)\\
	&\geq& -M_2(\kappa-1+\kappa o(1))
\end{eqnarray*}
for any constant $M_2>0$ and $\kappa>1.$ If we get $\liminf$ and limit of both sides of last inequality as $m, n\to\infty$ and $\kappa\to 1^+$ respectively, then we reach
\begin{equation*}
	\lim_{\kappa\to 1^+}\liminf_{m,n\to \infty}\min_{Q_n\leq Q_j\leq \kappa Q_{n}}\left(u_{mj}-u_{mn}\right)\geq0,
\end{equation*}
which means that $(u_{mn})$ is slowly decreasing relative to $(Q_n)$.  Therefore, we conclude with the help of Theorem \ref{maintheo-1} that $(u_{mn})$ is $P$-convergent to $\ell.$ \qed

\section{Main Results for the $(\overline{N}, p, q)$ Summable Double Sequences of complex numbers}\label{section4}

%\subsection{Main Results}\label{section3.1}\hspace*{\fill} \\

%\subsection{Examples}\label{section3.2}\hspace*{\fill} \\

%\subsection{Proofs of Main Results}\label{section3.3}\hspace*{\fill} \\
In parallel with the previous section, this section is also constructed by considering the following headings for $(\overline{N},p,q)$ summable double sequences of complex numbers:
\begin{itemize}
	\item[(a)] Determining certain subsets $\mathcal{T}\{u_{mn}\}$ of double sequence space and certain subclasses $\mathcal{C}\{p_m, q_n\}$ of single sequence space having the property: ``Under conditions  $\mathcal{C}\{p_m, q_n\},$ a double sequence $(u_{mn})\in \mathcal{T}\{u_{mn}\}$ being $(\overline{N},p,q)$ summable to $l\in \mathbb{R}$ is also $P$-convergent to the same value.''
	\item[(b)] Demonstrating with examples how the subsets $\mathcal{T}\{u_{mn}\}$ and the subclasses $\mathcal{C}\{p_m, q_n\}$ change for special means occurring depends on choosing of weights sequence $(p_{m})$ and $(q_{n}).$
	\item[(c)] Proving the property given in (a) for the related subsets and subclasses.
\end{itemize}
Here, we formulate our main results for double sequences of complex numbers as follows:
\begin{theorem}\label{maintheo-1*}
	Let $(p_m), (q_n)\in SV\hspace {-0.1 cm}A_{reg(\alpha)}$ and a double sequence $(u_{mn})$ be $(\overline{N},p,q)$ summable to a number $\ell.$
	If $(u_{mn})$ is slowly oscillating relative to $(P_m),$ slowly oscillating relative to $(Q_n)$, and slowly oscillating relative to $(P_m)$ or $(Q_n)$ in the strong sense, then $(u_{mn})$ is $P$-convergent to $\ell.$
\end{theorem}

\begin{theorem}\label{maintheo-2*}
	Let $(p_m), (q_n)\in SV\hspace {-0.1 cm}A_{reg(\alpha)}$ and a double sequence $(u_{mn})$ be $(\overline{N},p,q)$ summable to a number $\ell.$  If $(u_{mn})$ satisfies conditions
	\begin{equation}\label{hardytype(P)}
		\frac{P_{m}}{p_m}\Delta_{10}u_{mn}=O(1)\,\,\,\,\,\text{and}\,\,\,\,\,\frac{Q_{n}}{q_n}\Delta_{01}u_{mn}=O(1),
	\end{equation}
	then $(u_{mn})$ is $P$-convergent to $\ell.$	
\end{theorem}

\begin{exp-non}
	Some means occurring depends on choosing of weight sequences $(p_{m})$ and $(q_{m})$ are listed with their corresponding subsets $\mathcal{T}\{u_{mn}\}$ and subclasses $\mathcal{C}\{p_m, q_n\}$ for double sequences of complex numbers as follows.
	\begin{itemize}
		\item[(i)]  In case $p_{m}=q_n=1,$ it leads to the $(C,1,1)$ means of a double sequence where $P_{m}Q_n=(m+1)(n+1)$ for all $m, n\in\mathbb{N}.$ Under the circumstances, conditions in Theorem \ref{maintheo-1*} correspond to slow oscillation conditions in senses $(1,0), (0,1)$ and in the strong sense $(1,0)$ or $(0,1),$ besides conditions (\ref{hardytype(P)}) correspond to conditions
		\begin{equation*}
			m\Delta_{10}u_{mn}=O(1)\,\,\,\,\,\text{and}\,\,\,\,\,n\Delta_{01}u_{mn}=O(1),
		\end{equation*}
		however slow oscillation condition in sense $(1,1)$ and condition $mn\Delta_{11}$ $u_{mn}=O(1)$ discussed by M\'{o}ricz  \cite{moricz1994} therein are superfluous. As a result of our main Theorem \ref{maintheo-1*} and Theorem \ref{maintheo-2*}, the mentioned conditions are sufficient Tauberian conditions for the $(C,1,1)$ summability of double sequences of complex numbers.
		\item[(ii)] In case $\displaystyle{p_{m}q_n=\frac{1}{(m+1)(n+1)}},$ it leads to the \textit{logarithmic means}  of a double sequence where $P_{m}Q_n\sim\log m\log n$ for all $m, n\in\mathbb{N}.$ Under the circumstances, conditions in Theorem \ref{maintheo-1*} correspond to conditions of slow oscillation with respect to summability $(L,1)$ in sense $(1,0), (0,1)$ and the strong sense $(1,0)$ or $(0,1),$ i.e.,
		\begin{equation}\label{SOlog10}
			\lim_{\lambda\to 1^+}\limsup_{m,n\to \infty}\max_{m\leq i\leq m^\lambda}|u_{in}-u_{mn}|=0,
		\end{equation}
		\begin{equation}\label{SOlog01}
			\lim_{\kappa\to 1^+}\limsup_{m,n\to \infty}\max_{n\leq j\leq n^\kappa}|u_{mj}-u_{mn}|=0,
		\end{equation}
		and shall we say
		\begin{equation}\label{SOlog01strong}
			\lim_{\lambda, \kappa\to 1^+}\limsup_{m,n\to \infty}\max_{\substack{m\leq i\leq m^\lambda\\ n\leq j\leq n^\kappa}}|u_{ij}-u_{in}|= 0,
		\end{equation}
		which was discussed by Kwee \cite{kwee} and M\'{o}ricz \cite{moricz2013} for the logarithmic summability of single sequences in different ways. In addition, conditions (\ref{hardytype(P)}) correspond to conditions
		\begin{equation}\label{hardytypelog}
			m\log (m+1)\Delta_{10}u_{mn}=O(1)\,\,\,\,\,\text{and}\,\,\,\,\,n\log (n+1)\Delta_{01}u_{mn}=O(1).
		\end{equation}
		As a result of our main Theorem \ref{maintheo-1*} and Theorem \ref{maintheo-2*}, the mentioned conditions are sufficient Tauberian conditions for the logarithmic summability of double sequences of complex numbers.
	\end{itemize}
\end{exp-non}

Now, we present the proofs of our main results for double sequences of complex numbers.

\noindent \textbf{\emph{Proof of Theorem} \ref{maintheo-1*}}.
Assume that $(u_{mn})$ being $(\overline{N},p,q)$ summable to $\ell$ is slowly oscillating relative to $(P_m)$ and $(Q_n)$ and slowly oscillating relative to $(Q_n)$ in the strong sense. In order to prove that $(u_{mn})$ is $P$-convergent to same number, we indicate that the difference $|u_{mn}-\sigma_{mn}|$ is $P$-convergent to $0.$ By means of Lemma \ref{lemma}, we have for $\mu> m$ and $\eta>n$
\begin{eqnarray}\label{maintheo-eq1*}
	|u_{mn}-\sigma_{mn}|&\leq&\frac{P_{\mu}Q_{\eta}}{(P_{\mu}-P_m)(Q_{\eta}-Q_n)}\left|\sigma_{\mu\eta}-\sigma_{\mu n}-\sigma_{m\eta}+\sigma_{mn}\right|\nonumber\\
	&+&\frac{P_{\mu}}{P_{\mu}-P_m}\left|\sigma_{\mu n}-\sigma_{mn}\right|
	+\frac{Q_{\eta}}{Q_{\eta}-Q_n}\left|\sigma_{m\eta}-\sigma_{mn}\right|\nonumber\\
	&+&\frac{1}{(P_{\mu}-P_m)(Q_{\eta}-Q_n)}		 \sum_{i=m+1}^{\mu}\sum_{j=n+1}^{\eta}p_iq_j\left|u_{ij}-u_{mn}\right|\nonumber\\
	&\leq& \frac{P_{\mu}Q_{\eta}}{(P_{\mu}-P_m)(Q_{\eta}-Q_n)}\left|\sigma_{\mu\eta}-\sigma_{\mu n}-\sigma_{m\eta}+\sigma_{mn}\right|\nonumber\\
	&+&\frac{P_{\mu}}{P_{\mu}-P_m}\left|\sigma_{\mu n}-\sigma_{mn}\right|
	+\frac{Q_{\eta}}{Q_{\eta}-Q_n}\left|\sigma_{m\eta}-\sigma_{mn}\right|\nonumber\\
	&+&\max_{\substack{m\leq i\leq \mu\\ n\leq j\leq \eta}}	\left(u_{ij}-u_{in}\right)+\max_{m\leq i\leq\mu}\left(u_{in}-u_{mn}\right).
\end{eqnarray}
Put $\mu=\argmin\left\{P_i\geq\left(1+\displaystyle\frac{\delta}{2}\right) P_m \right\}$ and $\eta=\argmin\left\{Q_j\geq\left(1+\displaystyle\frac{\gamma}{2}\right)Q_n\right\}$ with $\delta, \gamma>0.$ Then, we get $\mu\geq m,\,\eta\geq n$ and $\displaystyle P_{\mu-1}<\left(1+\frac{\delta}{2}\right) P_m$ and $\displaystyle Q_{\eta-1}<\left(1+\frac{\gamma}{2}\right) Q_n.$ Moreover, on account of $(p_m), (q_n)\in SV\hspace {-0.1 cm}A_{reg(\alpha)},$ we have $\displaystyle \frac{P_m}{P_{m-1}}\leq 1+\frac{\delta}{4}$ and $\displaystyle \frac{Q_n}{Q_{n-1}}\leq 1+\frac{\gamma}{4}$ for sufficiently large $m, n.$ By means of Lemma \ref{lemma2}, it can be easily seen that $p_m/P_m\to 0$ and $q_n/Q_n\to 0$ as $m, n\to\infty,$ as well. As a result of these findings, we obtain by Lemma \ref{lemma2} that the inequalities (\ref{O(Pm)})-(\ref{O(Qn)*}) hold true for $\lambda,\,\kappa>1.$ From the assumption, we have $\sigma_{mn}\to\ell$ and so, $\sigma_{ij}-\sigma_{mn}\to 0$ for $i=m$ or $\mu$ and $j=n$ or $\eta$ as $m, n\to\infty.$ Since the sequences $(P_m)$ and $(Q_n)$ are strictly increasing sequences, we attain from (\ref{maintheo-eq1*}) that
\begin{eqnarray}\label{maintheo-eq2*}
	|u_{mn}-\sigma_{mn}| &\leq& \frac{P_{\mu}Q_{\eta}}{(P_{\mu}-P_m)(Q_{\eta}-Q_n)}\left|\sigma_{\mu\eta}-\sigma_{\mu n}-\sigma_{m\eta}+\sigma_{mn}\right|\nonumber\\
	&+&\frac{P_{\mu}}{P_{\mu}-P_m}\left|\sigma_{\mu n}-\sigma_{mn}\right|
	+\frac{Q_{\eta}}{Q_{\eta}-Q_n}\left|\sigma_{m\eta}-\sigma_{mn}\right|\nonumber\\&+&\max_{\substack{P_m\leq P_i\leq \lambda P_{m}\\ Q_n\leq Q_j\leq \kappa Q_{n}}}
	\left|u_{ij}-u_{in}\right|+\max_{P_m\leq P_i\leq \lambda P_{m}}\left|u_{in}-u_{mn}\right|\nonumber\\
	&\leq& \frac{4\lambda\kappa+2\lambda(\kappa-1)+2(\lambda-1)\kappa}{(\lambda-1)(\kappa-1)}(1+o(1)) o(1)\nonumber\\
	&+&\max_{\substack{P_m\leq P_i\leq \lambda P_{m}\\ Q_n\leq Q_j\leq \kappa Q_{n}}}
	\left|u_{ij}-u_{in}\right|+\max_{P_m\leq P_i\leq \lambda P_{m}}\left|u_{in}-u_{mn}\right|
\end{eqnarray}
where
\begin{multline*}
	\frac{P_\mu}{P_\mu-P_m}=\frac{P_\mu/P_m}{P_\mu/P_m-1}\leq\frac{2\lambda}{(\lambda-1)}(1+o(1))\,\,\,\,\,\,\,\,\text{and}\\
	\frac{Q_\eta}{Q_\eta-Q_n}=\frac{Q_\eta/Q_n}{Q_\eta/Q_n-1}\leq\frac{2\kappa}{(\kappa-1)}(1+o(1)).
\end{multline*}
If we get $\limsup$ of both sides of inequality (\ref{maintheo-eq2*}) as $m, n\to\infty,$ then we reach for any $\lambda, \kappa>1$
\begin{multline*}
	\limsup_{m, n\to\infty} |u_{mn}-\sigma_{mn}|\leq \limsup_{m, n\to\infty}\max_{\substack{P_m\leq P_i\leq \lambda P_{m}\\ Q_n\leq Q_j\leq \kappa Q_{n}}}
	\left|u_{ij}-u_{in}\right|\\+\limsup_{m, n\to\infty}\max_{P_m\leq P_i\leq \lambda P_{m}}\left|u_{in}-u_{mn}\right|.
\end{multline*}
If we get limit of both sides of last inequality as $\lambda, \kappa\to 1^+,$ then we find
\begin{equation}\label{maintheo-eq3*}
	\limsup_{m, n\to\infty} |u_{mn}-\sigma_{mn}|\leq 0
\end{equation}
due to the fact that $(u_{mn})$ is slowly oscillating relative to $(P_m)$ and $(Q_n)$ and slowly oscillating relative to $(Q_n)$ in the strong sense. Therefore, we conclude that $(u_{mn})$ is $P$-convergent to $\ell.$ \qed\\

\noindent \textbf{\emph{Proof of Theorem} \ref{maintheo-2*}}.
Assume that $(u_{mn})$ being $(\overline{N},p,q)$ summable to $\ell$ satisfies conditions (\ref{hardytype(P)}). If we indicate that conditions (\ref{hardytype(P)}) imply slow oscillation condition relative to $(P_m)$ and $(Q_n)$ and slow oscillation condition relative to $(Q_n)$ in the strong sense, then we could this proof with the help of Theorem \ref{maintheo-1*}. \\
Put $\mu=\argmin\left\{P_i\geq\left(1+\displaystyle\frac{\delta}{2}\right) P_m \right\}$ and $\eta=\argmin\left\{Q_j\geq\left(1+\displaystyle\frac{\gamma}{2}\right)  Q_n\right\}$ with $\delta, \gamma>0,$ so we can observe inequalities  (\ref{O(Pm)*})  and (\ref{O(Qn)*}). Then, we have for $n_0\leq m\leq i\leq\mu$ and $n_0\leq n$
\begin{eqnarray*}
	|u_{in}-u_{mn}|\leq\sum_{k=m+1}^{i}|\Delta_{10}u_{kn}|\leq M_1\sum_{k=m+1}^{i}\frac{p_k}{P_k}&\leq& M_1\left(\frac{P_\mu}{P_m}-1\right)\nonumber\\
	&\leq& M_1(\lambda-1+\lambda o(1))
\end{eqnarray*}
for any constant $M_1>0$ and $\lambda>1.$ If we get $\limsup$ and limit of both sides of last inequality as $m, n\to\infty$ and $\lambda\to 1^+$ respectively, then we reach
\begin{equation*}
	\lim_{\lambda\to 1^+}\limsup_{m,n\to \infty}\max_{P_m\leq P_i\leq \lambda P_{m}}\left|u_{in}-u_{mn}\right|=0,
\end{equation*}
which means that $(u_{mn})$ is slowly oscillating relative to $(P_m)$. Similarly, we obtain $n_0\leq n\leq j\leq\eta$ and $n_0\leq m$
\begin{eqnarray*}
	|u_{mj}-u_{mn}|\leq\sum_{r=n+1}^{j}|\Delta_{01}u_{mr}|\leq M_2\sum_{r=n+1}^{j}\frac{q_r}{Q_r}&\leq& M_2\left(\frac{Q_\eta}{Q_n}-1\right)\nonumber\\
	&\leq& M_2(\kappa-1+\kappa o(1))
\end{eqnarray*}
for any constant $M_2>0$ and $\kappa>1.$ If we get $\limsup$ and limit of both sides of last inequality as $m, n\to\infty$ and $\kappa\to 1^+$ respectively, then we reach
\begin{equation*}
	\lim_{\kappa\to 1^+}\limsup_{m,n\to \infty}\max_{Q_n\leq Q_j\leq \kappa Q_{n}}\left|u_{mj}-u_{mn}\right|=0,
\end{equation*}
which means that $(u_{mn})$ is slowly oscillating relative to $(Q_n)$.  Therefore, we conclude with the help of Theorem \ref{maintheo-1*} that $(u_{mn})$ is $P$-convergent to $\ell.$ \qed

\section*{Conclusion}

In this paper, we extended some theorems given for $(\overline{N},p)$ summable sequences to $(\overline{N},p,q)$ summable double sequences. We determined that conditions needed for $P-\lim\sigma_{mn}(u)=\ell$ to be $P-\lim u_{mn}=\ell$ are slow decrease (or oscillation) relative to $(P_m)$ and $(Q_n)$ and slow decrease (or oscillation) relative to $(P_m)$ or $(Q_n)$ in the strong sense under some additional condition imposed on $(p_m), (q_n).$ In the sequel, we presented a $O_L$-type (or $O$-type) Tauberian condition for $(\overline{N},p,q)$ summable sequences. Through this paper, we gave an answer to question (4) which was left as an open problem by Stadtm\"{u}ller and Baron in \cite{baronstadtmuller}. In addition to these, we indicate that our results include all of the classical Tauberian theorems for double sequences which $P$-convergence follows from Ces\`{a}ro and logarithmic summability under slow decrease (or oscillation) conditions relative to Schmidt and slow decrease (or oscillation) conditions relative to logarithmic summability in certain senses.

\section*{Acknowledgment}
%The authors are grateful to the anonymous referee for his/her valuable suggestions and comments in rewriting the article in the present form. The authors are also thankful to the Editor-in-Chief and the Associate Editor of the journal for their valuable comments.
The authors would like to thank Ulrich Stadtm\"{u}ller for sincerely answering our questions in reference to the problems which we encountered while conducting this study. The first author in this research is supported by The Scientific and Technological Research Council of Turkey (TUBITAK) under 2218 - National Postdoctoral Research Fellowship Program (Grant No. 118C577).

%%%%%%%%%%%%%%%%%%%%%%%%%%%%%%%%%%%%%%%%%

\end{document}